\documentclass[submission,copyright,creativecommons]{eptcs}
\providecommand{\event}{GASCOM 2026} 

\usepackage{amssymb}
\usepackage{amsmath}
\usepackage{cleveref}
\usepackage{amssymb}
\usepackage{amsmath}
\usepackage{amsthm}
\usepackage{graphicx}
\usepackage{enumerate}
\usepackage{tikz}
\usetikzlibrary{arrows}
\usetikzlibrary{arrows.meta}
\tikzset{>={Latex[width=1.2mm,length=1.7mm]}}
\usepackage{tabularx}
\usepackage{mathrsfs}

\usepackage{amssymb, latexsym}
\usepackage{graphicx}
\usepackage{xfrac}
\usepackage{verbatim}
\usepackage{mathtools}
\DeclarePairedDelimiter\ceil{\lceil}{\rceil}
\DeclarePairedDelimiter\floor{\lfloor}{\rfloor}

\DeclareSymbolFont{boldoperators}{OT1}{cmr}{bx}{n}
\SetSymbolFont{boldoperators}{bold}{OT1}{cmr}{bx}{n}

\newtheorem{thm}{Theorem}[section]
\newtheorem{prop}[thm]{Proposition}

\newtheorem{lem}[thm]{Lemma}
\newtheorem{obs}[thm]{Observation}
\newtheorem{conj}[thm]{Conjecture}
\newtheorem{prob}[thm]{Problem}
\theoremstyle{definition}
\newtheorem{defn}[thm]{Definition}
\newtheorem{alg}[thm]{Algorithm}

\numberwithin{equation}{section}

\newcommand{\sn}{\mathfrak{S}_n}

\newcommand{\A}{\mathcal{A}}
\newcommand{\D}{\mathcal{D}}
\newcommand{\U}{\mathcal{U}}

\newcommand{\imm}[1]{\mathrm{Imm}_{#1}}

\newcommand{\inv}{\textsc{inv}}

\newcommand{\defeq}{:=}

\newcommand{\sgn}{\mathrm{sgn}}

\newcommand{\ins}[1]{\mathrm{ins}_{#1}}
\newcommand{\inss}[1]{\mathrm{inss}_{#1}}

\newcommand{\tr}{{\negthickspace \top \negthickspace}}
\newcommand{\ntnsp}{\negthinspace}
\newcommand{\ntksp}{\negthickspace}

\newcommand{\bp}{\begin{prob}}
\newcommand{\ep}{\end{prob}}

\newcommand{\mat}[1]{\mathrm{Mat}_{#1 \times #1}}

\DeclareMathOperator{\perm}{per}

\newcommand{\permmon}[2]{#1_{1,#2_1} \ntnsp\cdots #1_{n,#2_n}}

\newcommand{\ssm}{\smallsetminus}

\newcommand{\net}[1]{\mathcal F_{#1}}

\newcommand{\vertex}[1]{#1}

\newcommand{\F}{\mathcal F}

\usepackage{iftex}

\ifpdf
  \usepackage{underscore}         
  \usepackage[T1]{fontenc}        
\else
  \usepackage{breakurl}           
\fi

\title{Permanental Inequalities and Unit Interval Orders}
\author{Sihong Pan \qquad Mark Skandera \qquad Jiayuan Wang
\institute{Lehigh University, Bethlehem PA, USA}
\email{sip318@lehigh.edu \quad mas906@lehigh.edu \quad jiw922@lehigh.edu}
}
\def\titlerunning{Permanental Inequalities and Unit Interval Orders}
\def\authorrunning{S. Pan, M. Skandera, \& J. Wang}
\begin{document}
\maketitle

\begin{abstract}
  Given an $n \times n$ matrix $A = (a_{i,j})$ and
  $I,J \subseteq [1,n] \defeq \{1,\dotsc,n\}$,
  let $A_{I,J} = (a_{i,j})_{i\in I, j \in J}$ denote the $(I, J)$-submatrix of $A$.
  We consider a class of $0$-$1$ totally nonnegative (TNN) matrices arising
  as antiadjacency matrices of unit interval orders, and show that for each
  matrix $A$ in this class,
  the inequalities
  \begin{equation}
    \perm(A_{[1,h],[1,h]}) \perm(A_{[h+1,n],[h+1,n]}) \;\geq\;
    \perm(A_{[1,n] \cap 2\mathbb Z, [1,n] \cap 2\mathbb Z})
    \perm(A_{[1,n] \ssm 2\mathbb Z, [1,n] \ssm 2\mathbb Z}) \tag{$*$} \label{eq:absmain}
  \end{equation}
  hold for $h = 1,\dotsc, n-1$. 

  Let $\mathfrak A_n$ be the Young Subgroup $\mathfrak{S}_{1, \dotsc, \lfloor\frac{n}{2}\rfloor}\times \mathfrak{S}_{\lceil\frac{n}{2}\rceil, \dotsc, n}$, and $\mathfrak B_n$ be the parity alternating permutations $\{w = w_1\cdots w_n \in \sn \ | \ \text{$i$ and $w_i$ have the same parity} \}$. We find a bijective map $f_n: \mathfrak A_n \rightarrow \mathfrak B_n$ such that $w \leq f_n(w)$ under the Bruhat order for $n\leq 13$ and conjecture this holds for all $n$. The existence of the bijective map $f_n$ ensures that inequalities~\eqref{eq:absmain} are satisfied for all TNN matrices $A$ when $h=\lfloor\frac{n}{2}\rfloor$. 
  We also conjecture the inequalities~\eqref{eq:absmain} to hold for all TNN matrices and all $h = 1, \dotsc ,n-1$. 
\end{abstract}

\section{Introduction}
Given an $n \times n$ matrix $A = (a_{i,j})$
and subsets $I,J \subseteq [n] \defeq \{1,\dotsc,n\}$,
define the submatrix $A_{I,J} = (a_{i,j})_{i \in I, j \in J}$ of $A$,
and call $\det(A_{I,J})$ a {\em minor} of $A$.
We call $A \in \mat n(\mathbb R)$ {\em totally nonnegative}
if each of its minors is negative.
It is known that the defining inequalities $\det(A_{I,J}) \geq 0$
of totally nonnegative matrices imply many other polynomial
inequalities in the matrix entries.
Some of these involve {\em immanants} of $A$ and its submatrices,
i.e., linear combinations of the $n!$ monomials
$\{ \permmon aw \,|\, w \in \sn \}$. 
In particular for any
function $f: \sn \rightarrow \mathbb R$, we define
the {\em $f$-immanant of $A$}~\cite{StanPos} to be the polynomial
\begin{equation}\label{eq:immdef}
  \imm f(A) \defeq \sum_{w \in \sn} f(w) \permmon aw
\end{equation}
generalizing the permanent and determinant
\begin{equation*}
  \perm(A) = \sum_{w \in \sn} \permmon aw,
  \qquad
  \det(A) = \sum_{w \in \sn} \sgn(w) \permmon aw.
  \end{equation*}

Skandera \cite{SkanIneq} gave a necessary and sufficient condition on subsets $I,I',J,J' \subseteq [n] := \{1,\dotsc, n\}$ such that the inequality 
\begin{equation*}
    \det(A_{I,J}) \det(A_{[n]\ssm I,[n] \ssm J}) \;\leq\;
    \det(A_{I',J'}) \det(A_{[n]\ssm I',[n] \ssm J'})
\end{equation*} holds whenever $A$ is TNN. Some other inequalities \cite{carlson1967,marcus1963perm,  StemImm, StemConj} suggest that one can use inequalities involving determinants to produce inequalities involving permanents by simply changing the direction of the inequality. 

Given a subset $I\subseteq [n]$, let $\bar I = [n] \ssm I$. A natural starting point for such an investigation
is Skandera's result~\cite[Thm.\,3.2]{SkanIneq}
which provides an {\em upper bound}
for all products $\det(A_{I,I}) \det(A_{\bar I, \bar I})$:
  \begin{equation*}
    \det(A_{I,I}) \det(A_{\bar I,\bar I}) \;\leq\;
    \det(A_{[n] \cap 2\mathbb Z, [n] \cap 2\mathbb Z})
    \det(A_{[n] \ssm 2\mathbb Z, [n] \ssm 2\mathbb Z})
  \end{equation*}
  whenever $A$ is TNN. 
  It is natural to ask whether we can find a {\em lower bound} for all
  products $\perm(A_{I,I}) \perm(A_{\bar I, \bar I})$:
\begin{equation}\label{eq:evenoddintro}
    \perm(A_{I,I}) \perm(A_{\bar I, \bar I}) \;\geq\;
    \perm(A_{[n] \cap 2\mathbb Z, [n] \cap 2\mathbb Z})
    \perm(A_{[n] \ssm 2\mathbb Z, [n] \ssm 2\mathbb Z}).
  \end{equation}

Our main result in Section~\ref{s:main} shows that inequalities of the form (\ref{eq:evenoddintro}) hold when
\begin{enumerate}
\item[(1)] The totally nonnegative matrix $A$ is related to a unit interval order, and
\item[(2)]
$I$ is the interval $[h]$ for any $h \in [n]$.
\end{enumerate}
We conjecture the inequalities
to hold for all totally nonnegative matrices and $I = [h]$.
On the other hand, we show that the inequalities do not hold for general $I$, even when restricted to matrices of the form (1).

\section{Four families of sets counted by Catalan numbers}\label{s:fourcatalan}

A simple subset of $n \times n$ totally nonnegative matrices is the set of
$0$-$1$
matrices which are weakly increasing in columns, weakly decreasing in rows,
and all entries on or below the diagonal are $1$.
Call this set $\A_n$.
For $A \in \A_n$, the boundary between entries of $A$ equal to $1$ and $0$
can be viewed as a $2n$-step {\em Dyck path},
consisting of $n$ horizontal
and $n$ vertical steps, none below the diagonal,
from the upper left corner of the
matrix to the lower right corner. 
For example, we can draw an $8$-step Dyck path on a matrix in $\A_4$. 

\begin{center}
\tikzset{every picture/.style={line width=0.75pt}} 

\begin{tikzpicture}[x=0.75pt,y=0.75pt,yscale=-.57,xscale=.57]

\draw  [dash pattern={on 0.84pt off 2.51pt}]  (280,80) -- (280,200) ;
\draw  [dash pattern={on 0.84pt off 2.51pt}][line width=0.75]  (120,40) -- (280,40) -- (280,200) -- (120,200) -- cycle ;
\draw [color={rgb, 255:red, 208; green, 2; blue, 27 }  ,draw opacity=1 ]   (120,40) -- (160,40) -- (160,80) -- (240,80) -- (240,160) -- (280,160) -- (280,200) ;

\draw (138,52.4) node [anchor=north west][inner sep=0.75pt]  [color={rgb, 255:red, 0; green, 0; blue, 0 }  ,opacity=1 ]  {$1$};
\draw (218,52.4) node [anchor=north west][inner sep=0.75pt]  [color={rgb, 255:red, 0; green, 0; blue, 0 }  ,opacity=1 ]  {$0$};
\draw (258,53.4) node [anchor=north west][inner sep=0.75pt]  [color={rgb, 255:red, 0; green, 0; blue, 0 }  ,opacity=1 ]  {$0$};
\draw (178,52.4) node [anchor=north west][inner sep=0.75pt]  [color={rgb, 255:red, 0; green, 0; blue, 0 }  ,opacity=1 ]  {$0$};
\draw (138,93.4) node [anchor=north west][inner sep=0.75pt]  [color={rgb, 255:red, 0; green, 0; blue, 0 }  ,opacity=1 ]  {$1$};
\draw (218,93.4) node [anchor=north west][inner sep=0.75pt]  [color={rgb, 255:red, 0; green, 0; blue, 0 }  ,opacity=1 ]  {$1$};
\draw (258,92.4) node [anchor=north west][inner sep=0.75pt]  [color={rgb, 255:red, 0; green, 0; blue, 0 }  ,opacity=1 ]  {$0$};
\draw (178,93.4) node [anchor=north west][inner sep=0.75pt]  [color={rgb, 255:red, 0; green, 0; blue, 0 }  ,opacity=1 ]  {$1$};
\draw (138,132.4) node [anchor=north west][inner sep=0.75pt]  [color={rgb, 255:red, 0; green, 0; blue, 0 }  ,opacity=1 ]  {$1$};
\draw (218,132.4) node [anchor=north west][inner sep=0.75pt]  [color={rgb, 255:red, 0; green, 0; blue, 0 }  ,opacity=1 ]  {$1$};
\draw (258,133.4) node [anchor=north west][inner sep=0.75pt]  [color={rgb, 255:red, 0; green, 0; blue, 0 }  ,opacity=1 ]  {$0$};
\draw (178,132.4) node [anchor=north west][inner sep=0.75pt]  [color={rgb, 255:red, 0; green, 0; blue, 0 }  ,opacity=1 ]  {$1$};
\draw (138,172.4) node [anchor=north west][inner sep=0.75pt]  [color={rgb, 255:red, 0; green, 0; blue, 0 }  ,opacity=1 ]  {$1$};
\draw (218,172.4) node [anchor=north west][inner sep=0.75pt]  [color={rgb, 255:red, 0; green, 0; blue, 0 }  ,opacity=1 ]  {$1$};
\draw (258,173.4) node [anchor=north west][inner sep=0.75pt]  [color={rgb, 255:red, 0; green, 0; blue, 0 }  ,opacity=1 ]  {$1$};
\draw (178,172.4) node [anchor=north west][inner sep=0.75pt]  [color={rgb, 255:red, 0; green, 0; blue, 0 }  ,opacity=1 ]  {$1$};

\end{tikzpicture}
\end{center}
The bijection between $\A_n$ and $2n$-step Dyck paths
implies that the cardinality of $\A_n$ is the $n$th Catalan number,
$\tfrac1{n+1}\tbinom{2n}{n}$.
(See \cite[Exercise 6.19 (i)]{StanEC2}.)

\begin{obs}\label{o:antiadjTNN}
  Each matrix in $\A_n$ 
  is totally nonnegative.
\end{obs}

Elements of $\A_n$ correspond bijectively to $n$-element posets
called {\em unit interval orders}.
Call this set $\U_n$.
(See \cite[p.\,33]{Fish}, \cite[Exercise 6.19 (ddd)]{StanEC2}, \cite[\S 8.2]{Trott}.)
These are the posets in which no induced four-element subposet
is isomorphic to either of the posets
\begin{equation*}
  \mathbf 3 + \mathbf 1 \;=\;
\begin{tikzpicture}[scale=.28,baseline=-3]
\draw[fill] (0,1) circle (1.2mm); 
\draw[fill] (0,0) circle (1.2mm); 
\draw[fill] (.8,0) circle (1.2mm); 
\draw[fill] (0,-1) circle (1.2mm); 
\draw[-] (0,-1) -- (0,0) -- (0,1);
\end{tikzpicture}\; ,
  \qquad
  \mathbf 2 + \mathbf 2 \;=\;
  \begin{tikzpicture}[scale=.50,baseline=-3]
\draw[fill] (0,.5) circle (1.2mm); 
\draw[fill] (0,-.5) circle (1.2mm); 
\draw[fill] (.8,-.5) circle (1.2mm); 
\draw[fill] (.8,.5) circle (1.2mm); 
\draw[-] (0,.5) -- (0,-.5);
\draw[-] (.8,.5) -- (.8,-.5);
\end{tikzpicture}\; .
\end{equation*}
The bijection $A \mapsto P(A)$ is given by
declaring $i \leq_{P(A)} j$ if and only if $a_{i,j} = 0$.
Each matrix $A$ is said to be the {\em antiadjacency matrix} of $P(A)$.

Define a {\em planar network of order $n$}
to be a directed, planar,
acyclic
multigraph
which can be embedded in a disc so that
$2n$ boundary vertices can be labeled clockwise as
{\em source} $1, \dotsc,$ {\em source} $n$,
{\em sink} $n, \dotsc,$ {\em sink} $1$.
We will assume all sources to have indegree $0$ and outdegree $1$,
and all sinks to have indegree $1$ and outdegree $0$.
Let $\net n$ denote the set of such networks.

For each subinterval $[a,b]$ of $[n]$ we define a {\em simple star network}
$F_{[a,b]} \in \net n$ by
\begin{enumerate}
\item Sources $1,\dotsc,n$ lie on a vertical line to the left;
  sinks $1,\dotsc,n$ lie on a vertical line to the right.
  Both are labeled from bottom to top.
\item An interior vertex lies between the sources and sinks.
\item For $i = 1, \dotsc, a-1$ and $i = b+1, \dotsc, n$, a directed edge
  begins at source $i$ and terminates at sink $i$.
\item For $i = a, \dotsc, b$, a directed edge begins at source $i$ and
  terminates at the interior vertex, and another directed edge begins
  at the interior vertex and terminates at sink $i$.
\end{enumerate}
For zero- and one-element subintervals
we define the trivial network $F_{\emptyset} = F_{[1,1]} = \cdots = F_{[n,n]}$
to have no interior vertex, and $n$ horizontal edges, each from source $i$
to sink $i$, for $i = 1,\dotsc,n$.
For example,
the element $F_{\smash{[2,4]}} \in \net 4$ is
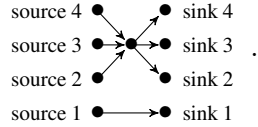
\begin{equation*}
\begin{tikzpicture}[scale=.45,baseline=15]
\node at (-2.5,2.5) {$\scriptstyle{ \mathrm{source}\;4}$};
\node at (-2.5,1.5) {$\scriptstyle{ \mathrm{source}\;3}$};
\node at (-2.5,0.5) {$\scriptstyle{ \mathrm{source}\;2}$};
\node at (-2.5,-0.5) {$\scriptstyle{ \mathrm{source}\;1}$};  
\node at (2.25,2.5) {$\scriptstyle{ \mathrm{sink}\;4}$};
\node at (2.25,1.5) {$\scriptstyle{ \mathrm{sink}\;3}$};
\node at (2.25,0.5) {$\scriptstyle{ \mathrm{sink}\;2}$};
\node at (2.25,-0.5) {$\scriptstyle{ \mathrm{sink}\;1}$};  
\draw[->,>=stealth'] (-1,2.5) -- (-.2,1.65);
\draw[->,>=stealth'] (-1,1.5) -- (-.2,1.5);
\draw[->,>=stealth'] (-1,0.5) -- (-.2,1.35);
\draw[->,>=stealth'] (0,1.5) -- (.8,2.35);
\draw[->,>=stealth'] (0,1.5) -- (.8,1.5);
\draw[->,>=stealth'] (0,1.5) -- (.8,0.65);
\draw[->,>=stealth'] (-1,-0.5) -- (.8,-0.5);
\node at (-1,2.5) {$\bullet$}; 
\node at (-1,1.5) {$\bullet$}; 
\node at (-1,0.5) {$\bullet$}; 
\node at (-1,-0.5) {$\bullet$};
\node at (1,2.5) {$\bullet$}; 
\node at (1,1.5) {$\bullet$}; 
\node at (1,0.5) {$\bullet$}; 
\node at (1,-0.5) {$\bullet$};
\node at (0,1.5) {$\bullet$};
\end{tikzpicture}
\ .
\end{equation*}  
For economy, we will omit edge orientations and vertices from
drawings of planar networks.
Thus we will draw the
seven simple star networks in $\net 4$ as
\begin{equation}\label{eq:simplestarnets}
\begin{gathered}
   \begin{tikzpicture}[scale=.3,baseline=-20]
\node at (-.4,0) {$\scriptstyle 4$};
\node at (-.4,-1) {$\scriptstyle 3$};
\node at (-.4,-2) {$\scriptstyle{2}$};
\node at (-.4,-3) {$\scriptstyle 1$};  
\node at (1.4,0) {$\scriptstyle 4$};
\node at (1.4,-1) {$\scriptstyle 3$};
\node at (1.4,-2) {$\scriptstyle{2}$};
\node at (1.4,-3) {$\scriptstyle 1$};  
\draw[-] (0,0) -- (1,-3);
\draw[-] (0,-1) -- (1,-2);
\draw[-] (0,-2) -- (1,-1);
\draw[-] (0,-3) -- (1,0);
\end{tikzpicture}
,\\
\phantom{\sum}\ntksp F_{[1,4]}\phantom{\sum}
\end{gathered}
\begin{gathered}
\begin{tikzpicture}[scale=.3,baseline=-20]
\node at (-.4,0) {$\scriptstyle 4$};
\node at (-.4,-1) {$\scriptstyle 3$};
\node at (-.4,-2) {$\scriptstyle{2}$};
\node at (-.4,-3) {$\scriptstyle 1$};  
\node at (1.4,0) {$\scriptstyle 4$};
\node at (1.4,-1) {$\scriptstyle 3$};
\node at (1.4,-2) {$\scriptstyle{2}$};
\node at (1.4,-3) {$\scriptstyle 1$};  
\draw[-] (0,0) -- (1,-2);
\draw[-] (0,-1) -- (1,-1);
\draw[-] (0,-2) -- (1,0);
\draw[-] (0,-3) -- (1,-3);
\end{tikzpicture},\\
   \phantom{\sum}\ntksp F_{[2,4]}\phantom{\sum}
 \end{gathered}
\begin{gathered}
\begin{tikzpicture}[scale=.3,baseline=-20]
\node at (-.4,0) {$\scriptstyle 4$};
\node at (-.4,-1) {$\scriptstyle 3$};
\node at (-.4,-2) {$\scriptstyle{2}$};
\node at (-.4,-3) {$\scriptstyle 1$};  
\node at (1.4,0) {$\scriptstyle 4$};
\node at (1.4,-1) {$\scriptstyle 3$};
\node at (1.4,-2) {$\scriptstyle{2}$};
\node at (1.4,-3) {$\scriptstyle 1$};  
\draw[-] (0,0) -- (1,0);
\draw[-] (0,-1) -- (1,-3);
\draw[-] (0,-2) -- (1,-2);
\draw[-] (0,-3) -- (1,-1);
\end{tikzpicture},\\
   \phantom{\sum}\ntksp F_{[1,3]}\phantom{\sum}
 \end{gathered}
\begin{gathered}
\begin{tikzpicture}[scale=.3,baseline=-20]
\node at (-.4,0) {$\scriptstyle 4$};
\node at (-.4,-1) {$\scriptstyle 3$};
\node at (-.4,-2) {$\scriptstyle{2}$};
\node at (-.4,-3) {$\scriptstyle 1$};  
\node at (1.4,0) {$\scriptstyle 4$};
\node at (1.4,-1) {$\scriptstyle 3$};
\node at (1.4,-2) {$\scriptstyle{2}$};
\node at (1.4,-3) {$\scriptstyle 1$};  
\draw[-] (0,0) -- (1,-1);
\draw[-] (0,-1) -- (1,0);
\draw[-] (0,-2) -- (1,-2);
\draw[-] (0,-3) -- (1,-3);
\end{tikzpicture},\\
   \phantom{\sum}\ntksp F_{[3,4]}\phantom{\sum}
 \end{gathered}
\begin{gathered}
\begin{tikzpicture}[scale=.3,baseline=-20]
\node at (-.4,0) {$\scriptstyle 4$};
\node at (-.4,-1) {$\scriptstyle 3$};
\node at (-.4,-2) {$\scriptstyle{2}$};
\node at (-.4,-3) {$\scriptstyle 1$};  
\node at (1.4,0) {$\scriptstyle 4$};
\node at (1.4,-1) {$\scriptstyle 3$};
\node at (1.4,-2) {$\scriptstyle{2}$};
\node at (1.4,-3) {$\scriptstyle 1$};  
\draw[-] (0,0) -- (1,0);
\draw[-] (0,-1) -- (1,-2);
\draw[-] (0,-2) -- (1,-1);
\draw[-] (0,-3) -- (1,-3);
\end{tikzpicture},\\
   \phantom{\sum}\ntksp F_{[2,3]}\phantom{\sum}
 \end{gathered}
\begin{gathered}
\begin{tikzpicture}[scale=.3,baseline=-20]
\node at (-.4,0) {$\scriptstyle 4$};
\node at (-.4,-1) {$\scriptstyle 3$};
\node at (-.4,-2) {$\scriptstyle{2}$};
\node at (-.4,-3) {$\scriptstyle 1$};  
\node at (1.4,0) {$\scriptstyle 4$};
\node at (1.4,-1) {$\scriptstyle 3$};
\node at (1.4,-2) {$\scriptstyle{2}$};
\node at (1.4,-3) {$\scriptstyle 1$};  
\draw[-] (0,0) -- (1,0);
\draw[-] (0,-1) -- (1,-1);
\draw[-] (0,-2) -- (1,-3);
\draw[-] (0,-3) -- (1,-2);
\end{tikzpicture},\\
   \phantom{\sum}\ntksp F_{[1,2]}\phantom{\sum}
 \end{gathered}
\
\begin{gathered}
\begin{tikzpicture}[scale=.3,baseline=-20]
\node at (-.4,0) {$\scriptstyle 4$};
\node at (-.4,-1) {$\scriptstyle 3$};
\node at (-.4,-2) {$\scriptstyle{2}$};
\node at (-.4,-3) {$\scriptstyle 1$};  
\node at (1.4,0) {$\scriptstyle 4$};
\node at (1.4,-1) {$\scriptstyle 3$};
\node at (1.4,-2) {$\scriptstyle{2}$};
\node at (1.4,-3) {$\scriptstyle 1$};  
\draw[-] (0,0) -- (1,0);
\draw[-] (0,-1) -- (1,-1);
\draw[-] (0,-2) -- (1,-2);
\draw[-] (0,-3) -- (1,-3);
\end{tikzpicture},\\
\phantom{\sum}\ntksp F_\emptyset \phantom{\sum}
 \end{gathered}
\end{equation}
where $F_\emptyset = F_{[1,1]} = F_{[2,2]} = F_{[3,3]} = F_{[4,4]}$.
When there is no danger of confusion, we will omit source and sink labels
as well. 

The set $\net n$ is generated by simple star networks and
two types of concatenation operations.
Given networks $E, F \in \net n$,
define
the (ordinary) concatenation $E \circ F$
as follows.  For $i = 1,\dotsc, n$, do
\begin{enumerate}
\item Remove sink $i$ of $E$ and source $i$ of $F$,
\item Merge each edge $(\vertex x, \text{sink }i)$ in $E$ with each edge
  $( \text{source }i, \vertex y )$ in $F$
  to form a single edge $(\vertex x, \vertex y)$ in $E \circ F$.
\end{enumerate}
Sometimes in a concatenation $E \circ F$,
there may exist vertices $\vertex x$ in $E$, $\vertex y$ in $F$
with more than one edge incident upon both.  Let
$m(\vertex x, \vertex y)$ be this number.
Define the {\em condensed concatenation}
$E \bullet F$ to be the
subdigraph of $E \circ F$ obtained
by removing, for all such pairs $(\vertex x, \vertex y)$,
all but one of the $m(\vertex x, \vertex y)$ edges incident upon both.
For example, in $\net 4$ we have 
\begin{equation}\label{eq:allbutone2}
  F_{[2,4]} \circ F_{[1,3]} =
\begin{tikzpicture}[scale=.35,baseline=0]
\node at (.6,1.5) {$\scriptstyle 4$};
\node at (.6,0.5) {$\scriptstyle 3$};
\node at (.6,-0.5) {$\scriptstyle{2}$};
\node at (.6,-1.5) {$\scriptstyle{1}$};  
\node at (3.4,1.5) {$\scriptstyle 4$};
\node at (3.4,0.5) {$\scriptstyle 3$};
\node at (3.4,-0.5) {$\scriptstyle{2}$};
\node at (3.4,-1.5) {$\scriptstyle{1}$};  
\draw[-] 
(1,1.5) -- (2,-0.5) -- (3,-0.5);
\draw[-] 
(1,-1.5) -- (2,-1.5) -- (3,0.5);
\draw[-] 
(1,-0.5) -- (2,1.5) -- (3,1.5);
\draw[-] 
(1,0.5) -- (2,0.5) -- (3,-1.5);
\end{tikzpicture}
\; \not \cong \;
F_{[2,4]} \bullet F_{[1,3]} =
\begin{tikzpicture}[scale=.35,baseline=0]
\node at (.6,1.5) {$\scriptstyle 4$};
\node at (.6,0.5) {$\scriptstyle 3$};
\node at (.6,-0.5) {$\scriptstyle{2}$};
\node at (.6,-1.5) {$\scriptstyle{1}$};  
\node at (3.4,1.5) {$\scriptstyle 4$};
\node at (3.4,0.5) {$\scriptstyle 3$};
\node at (3.4,-0.5) {$\scriptstyle{2}$};
\node at (3.4,-1.5) {$\scriptstyle{1}$};  
\draw[-] 
(1,1.5) -- (1.5,0.5) -- (2.5,-0.5) -- (3,-0.5);
\draw[-] 
(1,-1.5) -- (2,-1.5) -- (3,0.5);
\draw[-] 
(1,-0.5) -- (2,1.5) -- (3,1.5);
\draw[-] 
(1,0.5) -- (1.5,0.5) -- (2.5,-0.5) -- (3,-1.5);
\end{tikzpicture}\; .
\end{equation}


\begin{defn}\label{d:dsn}
  Define a {\em descending star network} to be a
  condensed concatenation of the form
  \begin{equation}\label{eq:dsn}
    F = F_{[c_1,d_1]} \bullet \cdots \bullet F_{[c_t,d_t]}, \\
\end{equation}
  with
    intervals $[c_1,d_1], \dotsc, [c_t,d_t]$
    distinct, pairwise nonnesting,
    and satisfying $c_1 > \cdots > c_t$.
\end{defn}
\noindent
Let $\D_n$ be the set of all descending star networks of order $n$.
In the notation of
Definition~\ref{d:dsn},
these are
$F_{[1,4]}$,
$F_{[2,4]} \bullet F_{[1,3]}$,
$F_{[2,4]} \bullet F_{[1,2]}$,
etc.
By \cite[Thm.\,3.6]{CHSSkanEKL}, \cite[Thm.\,3.5, Lem\,5.3]{SkanNNDCB},
the number of descending star networks of order $n$
is $\tfrac1{n+1} \tbinom{2n}n$.
Explicit bijections between
$\U_n$ and $\D_n$ can be found at~\cite{CHSSkanEhatKL}.  

Given a planar network $F$ of order $n$, define the {\em path matrix}
$B = (b_{i,j})$ of $F$
by
\begin{equation}\label{eq:pathmatrixdef}
  b_{i,j} = \text{number of paths in $F$ from source $i$ to sink $j$}.
\end{equation}
Let $\mathcal P_n$ be the set of $\tfrac 1{n+1}\tbinom{2n}n$ path matrices of
planar networks in $\mathcal D_n$.

Some, but not all, of the matrices in $\A_n$
appear in $\mathcal P_n$.
The precise connection between $\mathcal A_n$
and $\mathcal P_n$
can be described using \emph{block structures} of these matrices.
Let $\alpha = (\alpha_1,\dotsc,\alpha_r)$
be an integer composition of $n$, and define the (ordered) set partition
$(I_1, \dotsc, I_r)$ of $[n]$ 
 into intervals
 of lengths $\alpha_1,\dotsc,\alpha_r$ by
\begin{equation}\label{eq:alphaintervals}
  I_1 = [1,\alpha_1], \quad I_2 = [\alpha_1+1,\alpha_1+\alpha_2],\quad \dotsc, \quad
  I_r = [n-\alpha_r, n].
\end{equation} 
For any $n\times n$ matrix $A$, we can write it in terms of \emph{block submatrices} $A_{I_i, I_j}$ by
\small
\begin{equation}\label{eq:block}
  A = (A_{I_i,I_j}) =
  \begin{bmatrix}
    A_{I_1, I_1} & A_{I_1,I_2} & \cdots & A_{I_1,I_r} \\
    A_{I_2, I_1} & A_{I_2,I_2} & \cdots & A_{I_1,I_r} \\
    \vdots     & \vdots    &        & \vdots   \\
    A_{I_r, I_1} & A_{I_r,I_2} & \cdots & A_{I_r,I_r}
  \end{bmatrix}.
\end{equation}\normalsize
If $A_{I_i, I_j} =0$ whenever $i<j$, we say $A$ is
{\em $\alpha$-block lower-triangular};
if $A$ and $A^\tr$ are both
$\alpha$-block lower-triangular, we say that $A$ is
{\em $\alpha$-block diagonal}.
In particular, for any square matrix $A$ there is a
finest composition $\alpha = \alpha(A)$ of $n$ for which $A$ is $\alpha$-block
lower-triangular.
Next, we precisely state the relationship between $\mathcal A_n$
and $\mathcal P_n$.
Define the map
\begin{equation}\label{eq:psi}
  \begin{aligned}
    \psi: \mathcal A_n &\rightarrow \mathcal P_n\\
    A &\mapsto A_{I_1,I_1} \oplus \cdots \oplus A_{I_r,I_r},
  \end{aligned}
\end{equation}
where $(I_1,\dotsc,I_r)$ is the set partition associated 
to the finest composition $\alpha(A)$ as in (\ref{eq:alphaintervals}).

\begin{prop}\label{p:P A F psiA}
  Given unit interval order $P \in \mathcal U_n$
  with antiadjacency matrix $A$,  
  the descending star network $F = F(P) \in \mathcal D_n$
  has path matrix $\psi(A)$.
\end{prop}

\section{Permanents and path families}\label{s:permpathfam}

The planar networks defined in Section~\ref{s:fourcatalan} have close connections to permanents. 
Let $\pi=(\pi_1,\dotsc,\pi_n)$ be a sequence of source-to-sink paths in a
planar network $F \in \F_n$.
We call $\pi$ a \textit{path family} if there exists a permutation
$w=w_1\cdots w_n\in \sn$ such that $\pi_i$
is a path from source $i$ to sink $w_i$.
In this case, we say that $\pi$ has \emph{type} $w$.
Further, we say that the path family \emph{covers $F$}
if it contains every edge exactly once. 




Let $F$ be a planar network of order $n$ and let $I, J$ be subsets of $[n]$ with $|I|=|J|$. Let $\Pi_{I,J}(F)$ denote the set of path families from sources indexed by $I$ to sinks indexed by $J$. 
The following interpretation of
the determinent is
due to Karlin and McGregor~\cite{KMG}
and usually attributed to Lindstrom~\cite{LinVrep}.

\begin{lem}\label{c:llem}
  For $I, J \subseteq [n]$ with $|I| = |J| = k$, we have
  \begin{equation*}
    \det(A_{I,J}) = \# \{ \pi = (\pi_1, \dotsc, \pi_k) \in \Pi_{I,J}(F)
    \,|\, \pi_1,\dotsc,\pi_k \text{ pairwise nonintersecting } \}.
  \end{equation*}
\end{lem}  

We have the following combinatorial interpretation of the permanent using path families.

\begin{lem}\label{c:perlem}
    For $I, J \subseteq [n]$ with $|I| = |J| = k$, we have
    $\perm(A_{I,J}) = |\Pi_{I,J}(F)|$.
\end{lem}  

There is a simple formula counting path families covering a descending star network. 

\begin{lem}\label{eq:counting number of path families}
Given descending star network $F = F_{[c_1,d_1]} \bullet \cdots \bullet F_{[c_r,d_r]}$, define $m_1,\dotsc, m_r, k_1,\dotsc, k_{r-1}$ by
\begin{equation}\label{equ: size of stars}
        m_i=d_i-c_i+1, \quad
        k_i=\max\{d_{i+1}-c_i+1,0\}.
\end{equation}
Then we have 
\begin{equation}\label{equ: size}
    |\Pi_{[n],[n]}(F)| = \frac{m_1 !\cdots m_{r-1} !}{  k_{1} !\cdots k_{r-1}!} \cdot m_r!.
\end{equation}
\end{lem}

Let $h \in [n]$, then define sets $J = [h]$, $\bar J = [h+1,n]$, $I = [n] \cap 2 \mathbb Z$, and $\bar I = [n] \ssm 2\mathbb Z$. 
Lemma~\ref{eq:counting number of path families} and several factorial inequalities yield the following relation for path families.
\begin{lem}\label{cor: ascending star network}
    Fix $h\in[n]$ and let $F=F_{[c_1,d_1]} \bullet \cdots \bullet F_{[c_r,d_r]}$ be a descending star network such that $F_{[c_p,d_p]}$ is the only star with $[c_p,d_p] \cap [h,h+1] \neq \emptyset$.
    We have \begin{equation*}
        | \Pi_{I,I}(F)  | \cdot |\Pi_{\bar{I},\bar{I}}(F)| \leq | \Pi_{J,J}(F) | \cdot | \Pi_{\bar{J},\bar{J}}(F)  |.
    \end{equation*}
\end{lem}



\section{Merging stars in a star network}\label{s:merge}

For a general descending star network $F$, there may be multiple stars $F_{[c_p,d_p]}$ with $[c_p, d_p]\cap[h, h+1]\neq \emptyset$. We now define a function that can modify the star network $F$. For each number $h \in [n]$, define the \emph{$h$-merge function}
$\phi_h: \D_n \rightarrow \D_n$ as follows.
\begin{enumerate}
\setlength{\itemsep}{0.2em}
\item Given $F = F_{[c_1,d_1]} \bullet \cdots \bullet F_{[c_r,d_r]}$, define
  \begin{equation*}
    \begin{aligned}
      p &= \max \{k \,|\, [c_k,d_k] \cap [h,h+1] \neq \emptyset \},&
      q &= \min \{k \,|\, [c_k,d_k] \cap [h,h+1] \neq \emptyset \}.
    \end{aligned}
  \end{equation*}
\item Define $\phi_h(F) = F_{[c_1,d_1]} \bullet \cdots \bullet F_{[c_{q-1},d_{q-1}]} \bullet F_{[c_p,d_q]} \bullet F_{[c_{p+1},d_{p+1}]} \bullet \cdots \bullet F_{[c_r,d_r]}$.
\end{enumerate}

The $h$-merge function gives an equality on the cardinality of path families

\begin{prop}\label{p: star network contraction}
  Fix $h \in [n]$, and let $F$ be a descending star network of order $n$
  with path matrix $A$.
  Let $B$ be the path matrix of the descending star network
  $\phi_h(F)$. Then we have 
  \begin{enumerate}[(i)]
  \setlength{\itemsep}{0.7em}
  \item $A_{[h],[h]} = B_{[h],[h]}$ and
    $A_{[h+1,n],[h+1,n]} = B_{[h+1,n],[h+1,n]}$.
    
  \item $|\Pi_{[h],[h]}(F)|\cdot|\Pi_{[h+1,n],[h+1,n]}(F)| = |\Pi_{[h],[h]}(\phi_h(F))|\cdot|\Pi_{[h+1,n],[h+1,n]}(\phi_h(F))|$.
  \end{enumerate}
\end{prop}

\section{Main results}\label{s:main}

Now we state our main result about products of permanents.
\begin{thm}\label{t:main}
Let $A$ be the antiadjacency matrix of a unit interval order.
  Fix $h \in [n]$ and define sets
  \begin{equation*}
      J = [h] ,\qquad
      I = [n] \cap 2 \mathbb Z.
  \end{equation*}
  Then we have
  \begin{equation}\label{eq:main}
    \perm(A_{J,J})\perm(A_{\bar J, \bar J}) \geq \perm(A_{I,I})\perm(A_{\bar I, \bar I}).
  \end{equation}
  \end{thm}
  Proof omitted.

\begin{thm}\label{t:maintwo}
Let $A$ be the antiadjacency matrix of a unit interval order. Let $I = [n] \cap 2 \mathbb Z$. Then for any $i \in [n]$, we have 
    \begin{equation}\label{eq:secondmain}
        \perm(A_{[n]\setminus \{i\},[n]\setminus \{i\}})\perm(A_{\{i\},\{i\}}) \geq \perm(A_{I,I})\perm(A_{\bar I, \bar I}).
    \end{equation}
\end{thm}
Proof omitted.

\section{A Generalization on $I$}

A natural generalization of a product of two permanents indexed by even or odd integers is a product of $s$ permanents indexed by integers which are equivalent modulo $s$.
Define \begin{equation*}
    [c,d]_{r}^{s} = \{x\in [c,d] : \ x \equiv r \ (\text{mod } s) \}.  
\end{equation*}
Next, we have a generalization of Lemma~\ref{eq:counting number of path families}.
\begin{lem}\label{eq:counting r mod s} Let $F = F_{[c_0,d_0]} \bullet \cdots \bullet F_{[c_t,d_t]}$ be a descending star network. Fix $s\in[n]$. For every $r\in[0,s-1]$, We have
\begin{equation}\label{eq:evenoddfraction}
    |\Pi_{[1,n]_{r}^{s},[1,n]_{r}^{s}}(F)| = \frac{ |[c_0,d_0]_{r}^{s}|! \cdot |[c_1,d_1]_r^s|! \cdots |[c_t,d_t]_{r}^{s}|!}{|[c_0,d_1]_{r}^{s}|! \cdot |[c_1,d_2]_r^s|! \cdots |[c_{t-1},d_{t}]_r^s|!}.
\end{equation}
\end{lem}



We can now extend Theorem~\ref{t:main} to a more general result. 

\begin{thm}\label{t:rmodsmain}
   Let $A$ be the antiadjacency matrix of an $n$-element unit interval order. 
    
    
    Then for $s=2,\ldots, n-1$, we have 
    \begin{equation}\label{eq:new inequality}
   \perm(A_{[1,n]_0^s,[1,n]_0^s})\cdots \perm(A_{[1,n]_{s-1}^s,[1,n]_{s-1}^s})\geq \perm(A_{[1,n]_0^{s+1},[1,n]_0^{s+1}})\cdots\perm(A_{[1,n]_s^{s+1},[1,n]_s^{s+1}}).
    \end{equation}
    
\end{thm}
Proof omitted.


\section{A Generalization to all TNN matrices}

In this section, we will show that \eqref{eq:evenoddintro} holds for $I = [1,\lfloor \frac{n}{2} \rfloor]$. 

A polynomial $p(x_{11}, x_{12},\cdots, x_{nn})$ in $n^2$ variables is called totally non-negative (TNN) if it satisfies 
\begin{equation}\label{eq:TNNpolydef}
    p(A)\defeq p(a_{11}, a_{12},\cdots, a_{nn}) \geq 0
\end{equation} 
for each $n\times n$ totally non-negative matrix $A=(a_{ij})$. Therefore, \eqref{eq:evenoddintro} is true if and only if the polynomial \begin{equation}\label{TNNpoly}
    \perm(X_{[1,\lfloor \frac{n}{2} \rfloor],[1,\lfloor \frac{n}{2} \rfloor]}) \perm(X_{[\lceil \frac{n}{2} \rceil,n],[\lceil \frac{n}{2} \rceil,n]}) - \perm(X_{[n]\cap 2\mathbb Z , [n]\cap 2\mathbb Z})\perm(X_{[n]\setminus 2\mathbb Z , [n]\setminus 2\mathbb Z})
\end{equation} is TNN. 

Drake, Gerrish and Skandera~\cite{DGSBruhat} gave a necessary and sufficient condition for a polynomial which is a difference of ``permutation monomials" to be TNN. 

\begin{thm}[Drake, Gerrish, Skandera]\label{pairing test}
            Let $w$ and $u$ be two permutations in $\sn$. Then $w$ is less than or equal to $u$ in the Bruhat order if and only if the polynomial \[x_{1,w_1} \cdots x_{n,w_n} - x_{1,u_1} \cdots x_{n,u_n}\] is totally non-negative. 
\end{thm}

Let $t=\floor{\tfrac{n}{2}}$. The term $\perm(X_{\floor{\frac{n}{2}},\floor{\frac{n}{2}}})\perm(X_{\ceil{\frac{n}{2}},\ceil{\frac{n}{2}}})$ contains monomials 
$x_{1,w_1} \cdots x_{t,w_{t}} x_{t+1,w_{t+1}}\cdots x_{n,w_{n}}$ such that
\begin{equation*}
            w_1\cdots w_t \ \ \text{permutes} \ \{1,\dotsc,t\}, \qquad
            w_{t+1}\cdots w_n \ \ \text{permutes} \ \{t+1,\dotsc,n\}. 
        \end{equation*}
For convenience, we set 
\begin{equation}\label{def:An}
    \begin{split}
    &\mathfrak A_n = \mathfrak{S}_{1, \dotsc, t}\times \mathfrak{S}_{t+1, \dotsc, n} \\
    = \{&w \in \sn \ | \ w_1\cdots w_{t} \ \text{permutes} \ \{1,\dotsc,t\} \ \text{and} \ w_{t+1} \cdots w_n \ \text{permutes} \ \{t+1,\dotsc,n\} \}.
\end{split}
\end{equation}
Also, for future reference, we define \begin{equation}
    \tilde{\mathfrak A}_n = \{w \in \sn \ | \ w_1\cdots w_{t+1} \ \text{permutes} \ \{1,\dotsc,t+1\} \ \text{and} \ w_{t+2} \cdots w_n \ \text{permutes} \ \{t+2,\dotsc,n\} \}. 
\end{equation}
It is easy to see that $|\mathfrak A_n|= \floor{\tfrac{n}{2}}!\ceil{\tfrac{n}{2}}!$. 
For example, $\mathfrak A_4 = \{1234, 1243, 2134, 2143\}$. 
The product of permanents $\perm(X_{[n]\cap 2\mathbb{Z},[n]\cap 2\mathbb{Z}})\perm(X_{[n]\setminus 2\mathbb{Z}, [n]\setminus 2\mathbb{Z}})$ contains monomials    
$x_{1,u_1} \cdots x_{n,u_{n}}$ such that 
\begin{equation*}
            u_1u_3u_5\cdots  \ \text{permutes} \ \{1,3,5,\cdots\} = [n]\setminus 2\mathbb{Z}, \quad
            u_2u_4u_6\cdots  \ \text{permutes} \ \{2,4,6,\cdots\} = [n]\cap 2\mathbb{Z}.
        \end{equation*}       
We call the set that consists of $u = u_1\cdots u_n$ satisfying the above conditions $\mathfrak B_n$, i,e, \begin{equation}\label{def:Bn}
    \mathfrak B_n = \{u\in \sn \ | \ u_1u_2\cdots u_n \ \text{takes odd and even integers alternately and $u_1$ is odd} \}
\end{equation}
 It is easy to see that $|\mathfrak B_n|=\floor{\tfrac{n}{2}}! \ceil{\tfrac{n}{2}}!$.
Then, the polynomial~\eqref{TNNpoly} has this number 
of positive terms (on the left) and negative terms (on the right). Therefore, if we can find a bijective map $f_n: \mathfrak A_n \rightarrow \mathfrak B_n$, such that $w \leq f_n(w) $ for all $w\in \mathfrak A_n$, by Theorem~\ref{pairing test}, the polynomial~\eqref{TNNpoly} is TNN. 


Next, we introduce a few necessary maps. 
For $p, q\in[n]$, define an insertion map $\ins p$ and an insertion-swapping map $\inss q$ (when $n+1-q$ is even) given by
\begin{equation*}
\begin{aligned}
\ins p: \sn &\to \mathfrak{S}_{n+1},\\
w= w_1 \cdots w_p \cdots w_n&\mapsto w_1\cdots w_{p-1}(n+1) w_p \cdots w_n,\\
\text{inss}_q: \sn &\to \mathfrak{S}_{n+1},\\
w= w_1 \cdots w_q \cdots w_n&\mapsto w_1 \cdots w_{q-1} (n+1) w_{q+1}w_{q}w_{q+3}w_{q+2}\cdots w_{n}w_{n-1}.\\
\end{aligned}
\end{equation*}
For example, if $w=216435$, then $\ins 3(w) = 2176435$,  $\inss3(w) = 2174653$, but $\inss4(w)$ is not well-defined since $n+1-q=7-4=3$ is not even.

For $w\in \sn$, define the inversion statistic by
$\inv (w) \defeq |\{(w_i,w_j) \ | \ i < j \text{ and } w_i > w_j \}|$.
Observe that $|\mathfrak A_n|=|\tilde{\mathfrak A}_n|=\floor{\tfrac{n}{2}}!\ceil{\tfrac{n}{2}}!$. Then there exists a bijective map between $\mathfrak A_n$ and $\tilde{\mathfrak A}_n$. Let $w =  w_1 \cdots w_n \in \sn$. Define two maps $U$ and $R$ given by
\begin{equation*}
\begin{aligned}
U: \sn &\to \sn,\\
w &\mapsto w^U\defeq (n+1 -w_1)(n+1-w_2)\cdots (n+1-w_n),\\
R: \sn &\to \sn,\\
w &\mapsto w^R\defeq w_n \cdots w_1.\\
\end{aligned}
\end{equation*}
We abbreviate $U\circ R = R \circ U =RU$. This induces bijective maps between $\mathfrak A_n$ and $\tilde{\mathfrak A}_n$, and $\mathfrak B_n$ to itself.

 \begin{prop}\label{p: U and R}
 The map $U \circ R: \sn \rightarrow \sn$ sending $w \mapsto w^{UR}$ restricts to a bijection between
 $\mathfrak A_n$, $\tilde{\mathfrak A}_n$, and to an involution on $\mathfrak B_n$.
 \end{prop}

The map $U\circ R$ preserves the Bruhat order.

\begin{prop}\label{upside-down reverse}
    Let $v, \ w \in \sn$ with $v \leq w$. Then, $v^{RU} \leq w^{RU} $. 
\end{prop}





Next, we define a map $\tilde{f}_n: \tilde{\mathfrak A}_n \rightarrow \mathfrak B_n$. Let $\tilde{w} \in \tilde{\mathfrak A}_n$. Define $\tilde{f}_n$ by \begin{equation}\label{eq:def tilde f}
\begin{aligned}
    \tilde{f}_n : \tilde{\mathfrak A}_n &\rightarrow \mathfrak B_n, \\
    \tilde{w} &\mapsto [f_n(\tilde{w}^{RU})]^{RU}. 
    \end{aligned}
\end{equation}
We can show that $\tilde{f}_n$ is bijective if $f_n$ is bijective, and many properties of $f_n$ are inherited in $\tilde{f}_n$.



\begin{lem}\label{lem:induced bijection}
    Let $f_n: \mathfrak A_n \rightarrow \mathfrak B_n$ be a bijection such that $w \leq f_n(w)$ for all $w\in \mathfrak A_n$. Then, the bijection $\tilde{f}_n$~\eqref{eq:def tilde f} satisfies $\tilde{w} \leq \tilde{f}_n(\tilde{w})$ for every $\tilde{w} \in \tilde{\mathfrak A}_n$. 
\end{lem}

Now, we define $f_n:\mathfrak A_n \rightarrow \mathfrak B_n$ inductively. 

\begin{alg}\label{keyalg} 
Let $\mathfrak A_n$, $\tilde{\mathfrak A}_n$, $\mathfrak B_n$ be defined as above. 
and let $f_4$ be the map
\begin{equation}\label{alg base case}
    1234 \longmapsto 1234 \quad 1243 \longmapsto 1432 \quad 2134 \longmapsto 3214 \quad 2143 \longmapsto 3412.
\end{equation}

\noindent 
For $n=5,6,7,\dotsc,$ do
\begin{enumerate}
    \item
If $n$ is odd, do 
\begin{enumerate}[a.]
    \item Write $n-1 = 2k$.
    \item For each $w\in \mathfrak A_n$, find the unique pair $(a,p) \in \mathfrak A_{n-1} \times \mathbb N$ such that $w = \ins p (a)$. 
    \item Define 
    $f_n(w) = \inss{2(p-k)-1}(f_{n-1}(a))$.
\end{enumerate}
\item If $n$ is even, do 
\begin{enumerate}[a.]
    \item Write $n-1 = 2k+1$. 
    \item Define $\tilde{f}_{n-1} : \tilde{\mathfrak A}_n \rightarrow \mathfrak B_n$ by $f_{n-1}$ as in~\eqref{eq:def tilde f}. 
    \item For each $u \in \mathfrak A_{n}$, $(\tilde{w},q) \in \tilde{\mathfrak A}_n \times \mathbb N$ such that $u=\ins q (\tilde{w})$. 
    \item[d.] Define
      $f_n(u) = \inss{2(q-k-1)}(\tilde{f}_{n-1}(\tilde{w}))$.
\end{enumerate}
\end{enumerate}
\end{alg}

\noindent
We illustrate the inductive definition of $f_n$ in Algorithm~\ref{keyalg} with the following diagram.


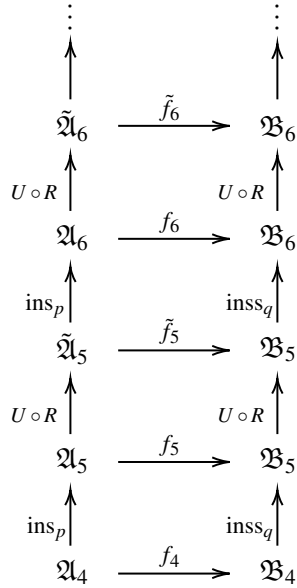
\begin{center}

\tikzset{every picture/.style={line width=0.75pt}} 

\begin{tikzpicture}[x=0.75pt,y=0.75pt,yscale=-0.7,xscale=0.8]

\draw    (300,380) -- (300,342) ;
\draw [shift={(300,340)}, rotate = 90] [color={rgb, 255:red, 0; green, 0; blue, 0 }  ][line width=0.75]    (10.93,-3.29) .. controls (6.95,-1.4) and (3.31,-0.3) .. (0,0) .. controls (3.31,0.3) and (6.95,1.4) .. (10.93,3.29)   ;
\draw    (170,380) -- (170,342) ;
\draw [shift={(170,340)}, rotate = 90] [color={rgb, 255:red, 0; green, 0; blue, 0 }  ][line width=0.75]    (10.93,-3.29) .. controls (6.95,-1.4) and (3.31,-0.3) .. (0,0) .. controls (3.31,0.3) and (6.95,1.4) .. (10.93,3.29)   ;
\draw    (300,460) -- (300,422) ;
\draw [shift={(300,420)}, rotate = 90] [color={rgb, 255:red, 0; green, 0; blue, 0 }  ][line width=0.75]    (10.93,-3.29) .. controls (6.95,-1.4) and (3.31,-0.3) .. (0,0) .. controls (3.31,0.3) and (6.95,1.4) .. (10.93,3.29)   ;
\draw    (170,460) -- (170,422) ;
\draw [shift={(170,420)}, rotate = 90] [color={rgb, 255:red, 0; green, 0; blue, 0 }  ][line width=0.75]    (10.93,-3.29) .. controls (6.95,-1.4) and (3.31,-0.3) .. (0,0) .. controls (3.31,0.3) and (6.95,1.4) .. (10.93,3.29)   ;
\draw    (170,140) -- (170,101) ;
\draw [shift={(170,99)}, rotate = 90] [color={rgb, 255:red, 0; green, 0; blue, 0 }  ][line width=0.75]    (10.93,-3.29) .. controls (6.95,-1.4) and (3.31,-0.3) .. (0,0) .. controls (3.31,0.3) and (6.95,1.4) .. (10.93,3.29)   ;
\draw    (300,140) -- (300,102) ;
\draw [shift={(300,100)}, rotate = 90] [color={rgb, 255:red, 0; green, 0; blue, 0 }  ][line width=0.75]    (10.93,-3.29) .. controls (6.95,-1.4) and (3.31,-0.3) .. (0,0) .. controls (3.31,0.3) and (6.95,1.4) .. (10.93,3.29)   ;
\draw    (200,480) -- (268,480) ;
\draw [shift={(270,480)}, rotate = 180] [color={rgb, 255:red, 0; green, 0; blue, 0 }  ][line width=0.75]    (10.93,-3.29) .. controls (6.95,-1.4) and (3.31,-0.3) .. (0,0) .. controls (3.31,0.3) and (6.95,1.4) .. (10.93,3.29)   ;
\draw    (200,400) -- (268,400) ;
\draw [shift={(270,400)}, rotate = 180] [color={rgb, 255:red, 0; green, 0; blue, 0 }  ][line width=0.75]    (10.93,-3.29) .. controls (6.95,-1.4) and (3.31,-0.3) .. (0,0) .. controls (3.31,0.3) and (6.95,1.4) .. (10.93,3.29)   ;
\draw    (200,159) -- (268,159) ;
\draw [shift={(270,159)}, rotate = 180] [color={rgb, 255:red, 0; green, 0; blue, 0 }  ][line width=0.75]    (10.93,-3.29) .. controls (6.95,-1.4) and (3.31,-0.3) .. (0,0) .. controls (3.31,0.3) and (6.95,1.4) .. (10.93,3.29)   ;
\draw    (300,220) -- (300,182) ;
\draw [shift={(300,180)}, rotate = 90] [color={rgb, 255:red, 0; green, 0; blue, 0 }  ][line width=0.75]    (10.93,-3.29) .. controls (6.95,-1.4) and (3.31,-0.3) .. (0,0) .. controls (3.31,0.3) and (6.95,1.4) .. (10.93,3.29)   ;
\draw    (170,220) -- (170,182) ;
\draw [shift={(170,180)}, rotate = 90] [color={rgb, 255:red, 0; green, 0; blue, 0 }  ][line width=0.75]    (10.93,-3.29) .. controls (6.95,-1.4) and (3.31,-0.3) .. (0,0) .. controls (3.31,0.3) and (6.95,1.4) .. (10.93,3.29)   ;
\draw    (300,300) -- (300,262) ;
\draw [shift={(300,260)}, rotate = 90] [color={rgb, 255:red, 0; green, 0; blue, 0 }  ][line width=0.75]    (10.93,-3.29) .. controls (6.95,-1.4) and (3.31,-0.3) .. (0,0) .. controls (3.31,0.3) and (6.95,1.4) .. (10.93,3.29)   ;
\draw    (170,300) -- (170,262) ;
\draw [shift={(170,260)}, rotate = 90] [color={rgb, 255:red, 0; green, 0; blue, 0 }  ][line width=0.75]    (10.93,-3.29) .. controls (6.95,-1.4) and (3.31,-0.3) .. (0,0) .. controls (3.31,0.3) and (6.95,1.4) .. (10.93,3.29)   ;
\draw    (200,320) -- (268,320) ;
\draw [shift={(270,320)}, rotate = 180] [color={rgb, 255:red, 0; green, 0; blue, 0 }  ][line width=0.75]    (10.93,-3.29) .. controls (6.95,-1.4) and (3.31,-0.3) .. (0,0) .. controls (3.31,0.3) and (6.95,1.4) .. (10.93,3.29)   ;
\draw    (200,240) -- (268,240) ;
\draw [shift={(270,240)}, rotate = 180] [color={rgb, 255:red, 0; green, 0; blue, 0 }  ][line width=0.75]    (10.93,-3.29) .. controls (6.95,-1.4) and (3.31,-0.3) .. (0,0) .. controls (3.31,0.3) and (6.95,1.4) .. (10.93,3.29)   ;

\draw (171.5,479) node    {$\mathfrak A_{4}$};
\draw (300.5,479) node    {$\mathfrak B_{4}$};
\draw (300.5,399) node    {$\mathfrak B_{5}$};
\draw (171.5,399) node    {$\mathfrak A_{5}$};
\draw (175.15,68) node [anchor=north west][inner sep=0.75pt]  [rotate=-90]  {$\cdots $};
\draw (305.15,68) node [anchor=north west][inner sep=0.75pt]  [rotate=-90]  {$\cdots $};
\draw (141,437.4) node [anchor=north west][inner sep=0.75pt]  [font=\footnotesize]  {$\ins p$};
\draw (266,437.4) node [anchor=north west][inner sep=0.75pt]  [font=\footnotesize]  {$\inss q$};
\draw (131,359.4) node [anchor=north west][inner sep=0.75pt]  [font=\scriptsize]  {$\mathnormal{{\textstyle U\circ R}}$};
\draw (261,359.4) node [anchor=north west][inner sep=0.75pt]  [font=\scriptsize]  {$\mathnormal{{\textstyle U\circ R}}$};
\draw (224,458.4) node [anchor=north west][inner sep=0.75pt]  [font=\footnotesize]  {$f_{4}$};
\draw (224,378.4) node [anchor=north west][inner sep=0.75pt]  [font=\footnotesize]  {$f_{5}$};
\draw (171.5,319) node    {$\tilde{\mathfrak A}_{5}$};
\draw (300.5,319) node    {$\mathfrak B_{5}$};
\draw (300.5,239) node    {$\mathfrak B_{6}$};
\draw (171.5,239) node    {$\mathfrak A_{6}$};
\draw (140,277.4) node [anchor=north west][inner sep=0.75pt]  [font=\footnotesize]  {$\ins p$};
\draw (266,277.4) node [anchor=north west][inner sep=0.75pt]  [font=\footnotesize]  {$\inss q$};
\draw (131,199.4) node [anchor=north west][inner sep=0.75pt]  [font=\scriptsize]  {$\mathnormal{{\textstyle U\circ R}}$};
\draw (261,199.4) node [anchor=north west][inner sep=0.75pt]  [font=\scriptsize]  {$\mathnormal{{\textstyle U\circ R}}$};
\draw (224,295) node [anchor=north west][inner sep=0.75pt]  [font=\footnotesize]  {$\tilde{f}_{5}$};
\draw (224,218.4) node [anchor=north west][inner sep=0.75pt]  [font=\footnotesize]  {$f_{6}$};
\draw (171.5,159) node    {$\tilde{\mathfrak A}_{6}$};
\draw (300.5,161) node    {$\mathfrak B_{6}$};
\draw (224,135) node [anchor=north west][inner sep=0.75pt]  [font=\footnotesize]  {$\tilde{f}_{6}$};

\end{tikzpicture}

\end{center}

Our algorithm behaves nicely with respect to inversion and the Bruhat order. We have the following partial results and conjecture for the function $f_n$ defined by Algorithm \ref{keyalg}. 
\begin{thm}\label{l:nec Bruhat condition}
    For all $w\in \mathfrak A_n \setminus \{e\}$, we have $\inv (w) < \inv (f_n(w))$.
\end{thm}
Proof omitted.

\begin{thm}\label{t:special n}
    For $n\leq 13$ and all $w 
    \in \mathfrak A_n$, we have $w\leq f_n(w)$ in the Bruhat order.
\end{thm}
Proof omitted.

\begin{conj}\label{conj:general}
    For all $n$ and all $w \in \mathfrak A_n$ we have $w \leq f_n(w)$ in the Bruhat order.
\end{conj}

A proof of Conjecture~\ref{conj:general} 
would extend Theorem~\ref{t:main} to all totally nonnegative matrices.

\bibliographystyle{eptcs}
\bibliography{skan}
\end{document}